\documentclass{article}

\usepackage{amsmath}
\usepackage{amsfonts}

\newtheorem{theo}{Theorem}

\newtheorem{prop}{Proposition}
\newtheorem{lemm}{Lemma}

\def\N{\mathbb{N}}          

\def\0{{\bf 0}}

\def\R{\mathbb{R}}

\def\Ed{{\cal E}}

\def\essinf{{\rm ess~inf}}

\newcommand{\tod}{\stackrel{{\cal D}}{\longrightarrow}}

\newcommand{\eqco}{\setcounter{equation}{0}}
\newcommand{\thco}{\setcounter{theo}{0}}
\newcommand{\prco}{\setcounter{prop}{0}}
\newcommand{\laco}{\setcounter{lemm}{0}}
\newcommand{\coco}{\setcounter{coro}{0}}
\newcommand{\cjco}{\setcounter{conj}{0}}

\newcommand{\deco}{\setcounter{defn}{0}}
\newcommand{\allco}{\eqco  \thco \prco \laco \coco \cjco \deco}
\newcommand{\qed}{\rule[-1mm]{3mm}{3mm}}

\newcommand{\Po}{{\cal P}}
\newcommand{\Q}{{\cal Q}}
\renewcommand{\NG}{{\rm NG}}
\newcommand{\MST}{{\rm MST}}
\newcommand{\Vor}{{\rm VOR}}
\newcommand{\SIG}{{\rm SIG}}
\newcommand{\Del}{{\rm DEL}}

\renewcommand{\H}{{\cal H}}
\renewcommand{\P}{{\cal P}}
\renewcommand{\SS}{{\cal S}}
\newcommand{\A}{{\cal A}}

\newcommand{\essup}{{\rm ess~sup}}

\newcommand{\limn}{\lim_{n \to \infty} }

\newcommand{\X}{{\cal X}}
\newcommand{\K}{{\cal K}}
\newcommand{\F}{{\cal F}}

\newcommand{\eps}{\varepsilon}
\def\bdm{\begin{displaymath}}
\newcommand{\edm}{\end{displaymath}}
\def\benu{\begin{enumerate}}
\def\eenu{\end{enumerate}}
\def\beqn{\begin{equation}}
\def\eeqn{\end{equation}}
\def\be{\begin{equation}}
\def\ee{\end{equation}}
\def\bea{\begin{eqnarray}}
\def\eea{\end{eqnarray}}
\newcommand{\bean}{\begin{eqnarray*}}
\newcommand{\eean}{\end{eqnarray*}}
\newcommand{\bear}{\begin{eqnarray}}
\newcommand{\eear}{\end{eqnarray}}

\renewcommand{\epsilon}{\varepsilon}

\begin{document}

\title{\bf Weak laws in geometric probability}

\bigskip

\author{Mathew D. Penrose and J. E. Yukich$^{1}$ \\
\\
{\normalsize{\em University of Durham and Lehigh University}} }

\date{July 2001}
\maketitle

\footnotetext[1] { Department of Mathematical Sciences, University
of Durham, South Road, Durham DH1 3LE, England: {\texttt
mathew.penrose@durham.ac.uk} }

\footnotetext[2]{ Department of Mathematics, Lehigh University,
Bethlehem PA 18015, USA: {\texttt joseph.yukich@lehigh.edu} }

\footnotetext{$~^1$ Research supported in part by NSA grant
MDA904-01-1-0029 }

\begin{abstract} Using a  coupling argument, we  establish
 a general weak law of large numbers
for functionals of binomial
  point processes in $d$-dimensional space, with a limit that
depends explicitly on the (possibly non-uniform) density of the
point process.  The general result is applied to the minimal
spanning tree, the $k$-nearest neighbors graph,  the Voronoi
graph, and the sphere of influence graph.  Functionals of
interest include total edge length with arbitrary weighting,
number of vertices of specifed degree, and number of components.
We also obtain weak laws for functionals of marked point
processes, including statistics of Boolean models.

\end{abstract}

\section{ Introduction}

Establishing laws of large numbers (LLN) for functionals of
random Euclidean point sets is of considerable interest.  When
the point set forms the vertex set of a graph, functionals of
interest include total edge length with arbitrary weighting,
number of edges, and number of components.  Relevant graphs
include those in computational geometry, such as the minimal
spanning tree, the $k$-nearest neighbors graph, Voronoi graph,
and sphere of influence graph (these graphs are formally defined
in Section \ref{secapp}). When the random Euclidean point set is
a marked point set, then functionals of interest include those
arising in the stochastic geometry of Boolean models.

For many functionals, subadditivity works well as a basic tool;
see \cite{Stbk} and \cite{Y1}  for surveys. This is the case with
power-weighted edge length functionals. For example, if $G(n,f)$
is the minimal spanning tree (MST)
 on $n $ i.i.d. random
$d$-vectors
 with common
density $f$ on $\R^d$, and if $|e|$ denotes the length of the
edge $e$, then  the following asymptotics hold for the sum
$\sum_{e\in G(n,f)} \vert e \vert^\alpha$ of the power-weighted
edge lengths:

\begin{theo}\label{basicMST}  If either (i) $1 \leq \alpha < d,
\int_{\R^d} f(x)^{(d-\alpha)/d} dx < \infty$,  and $\int_{\R^d}
\vert x \vert^r f(x) dx < \infty$ for some $r
> d/(d-\alpha)$
 or (ii) $\alpha \geq d$ and  $f$ has support on $[0,1]^d$,
 and is bounded away from zero,  then as $n \to \infty$
   \bea  n^{(\alpha-d)/d}
\sum_{e \in G(n,f)} \vert e \vert^\alpha \to C(\alpha,d) \int_{\R^d}
f(x)^{(d-\alpha)/d} dx \ \ a.s., \label{bhh} \eea where $C(\alpha,d)$ is a
positive constant.
\end{theo}
For a proof see  \cite{Y1} for case (i) and see \cite{Y2} for
case (ii).

However, many functionals are not amenable to subadditive methods.
For example, if the edge lengths $\vert e \vert$ in a Euclidean
graph $G$ are weighted by a general function $\phi$, giving a sum
of the form $\sum_{e \in G} \phi(\vert e \vert)$, then
subadditive methods break down, and the LLN behavior is much less
well understood. Functionals involving the Voronoi, Delaunay, and
sphere of influence graphs are also generally not amenable to
subadditive methods.  Functionals of marked point processes are
usually not subadditive either. Moreover, even when subadditive
methods are applicable, they provide little information on the
numerical values of limiting constants such as $C(\alpha,d)$
appearing in (\ref{bhh}).

An alternative approach is the so-called ``objective method''.
Steele \cite{Stbk} coined this term for
a  
 philosophy whereby,
 loosely speaking, one uses the locally Poisson nature of a binomial
 point process to  describe the limiting behavior of functionals on finite
 point sets in terms of
 related functionals  defined on infinite Poisson point sets.
 Aldous and Steele \cite{AS} used this idea to analyze certain functionals
associated with the MST on uniform points, but one might expect
it to be applicable to any functional, including those defined
over non-uniform point samples,  consisting of contributions which
are locally determined in some sense. As noted in \cite{AS,Stbk},
 making formal sense of this intuition is not always quite
so simple as one might imagine.

In an attempt to formulate in general terms the idea of locally determined
contributions,
Penrose and Yukich \cite{PY1,PY2} introduce a concept
 of ``stabilizing'' functionals
and essentially use the objective method to establish a strong
law of large numbers (Theorem 3.2 of \cite{PY2}) for stabilizing
functionals on uniform point sets in $\R^d$.  Jimenez and Yukich
\cite{JiY}  obtain sufficient conditions yielding laws of large
numbers for sums involving general edge weights and non-uniform
point sets,
 but their conditions are rather strong and are limited to
 increasing functions $\phi$.

The goal of this paper is to use the objective method to provide
relatively simple conditions guaranteeing a general weak LLN  for
stabilizing functionals on possibly non-uniform point samples of
size $n$. We illustrate the diverse applications of the general
LLN by obtaining weak laws
 for functionals of spatial point processes  in
computational geometry as well as functionals of marked point
processes, including those arising in packing processes and the
stochastic geometry of Boolean models.

 Many stabilizing functionals are defined in terms of
   graphs $G$ which are themselves stabilizing, i.e., locally determined
in a sense to be made precise below; stabilizing graphs include
 the  MST, $k$-nearest neighbor,
Voronoi, and sphere of influence graph.
 Given  a stabilizing graph $G$,
the theory applies to functionals such as
 the number of leaves, the number
 of components, and the
 sum of weighted edge lengths $\sum_{e \in G}
\phi(\vert e \vert)$.
In these graphs, edges are between
 ``nearby" points, and since the density of points
grows in proportion to $n$, the typical distance between nearby
points can be thought of as decreasing in proportion to
$n^{-1/d}$. Therefore we
 consider sums of the form $\sum_{e \in G}\phi(n^{1/d}|e|)$,
establishing weak LLN behavior.
 The limiting constants are defined
explicitly in terms of $\phi$,  the density $f$, and certain  graphs
on Poisson processes, thereby providing extra information
on the value of limiting  constants such as that arising in (\ref{bhh})
that is not given by subadditive methods alone, even in
the classic case when
 $\phi$ is the identity function.

\section{Main Results}
\label{sec1} \allco

\subsection{Terminology}
\label{secterm}

In Section \ref{gensec} we shall  formulate
 a collection of general LLN results.
 Before doing so we need some terminology.

Given $\X \subset \R^d$ and a positive scalar $a$, let $a\X:=
\{ax: x \in \X\}.$ Given $y \in \R^d$ set $y+\X :=\{y+x: x \in
\X\}$. For $x \in \R^d$,
let $|x|$ be its Euclidean modulus
and for
 $r>0$, let $B(x;r)$ denote the
Euclidean ball $\{ y \in \R^d: |y-x| \leq r\}$.
Let $\0$ denote the origin of $\R^d$.

Suppose $\xi(x; \X)$ is a  measurable
   $\R^+$-valued function
defined for all pairs $(x,\X)$, where $\X \subset \R^d$  and $x$
is an element of $\X$. We assume that $\X$ is locally finite,
i.e., contains only finitely many points in any bounded region.
 Suppose $\xi$ is translation invariant,
i.e. $\xi(y+x;y+\X)= \xi(y;\X)$ for all $y \in \R^d$ and all
$x,\X$.
 Then $\xi$   induces a
translation-invariant functional $H_{\xi}$ defined on finite
point sets $\X \subset \R^d$ by
  \bea H_{\xi}(\X) :=
\sum_{x \in \X} \xi(x; \X). \label{induceh} \eea Functionals
admitting the representation (\ref{induceh})
 include the total edge length, the total number of edges,  the total
 number of components, and total number of vertices of fixed degree
of Euclidean graphs. Later on, in cases with $x \notin \X$ it
will be useful to abbreviate the notation $\xi(x;\X \cup \{x\})$
to
 $\xi( x; \X)$.


We probe the behavior of the  functional $H_{\xi}$
 by inserting an  extra point into its domain.
For ``typical"  point sets $\X$, it is conceivable that the
contribution $\xi(x; \X)$ is not affected by changes in
 $\X$ which are far
from $x$.
 We formalize this notion 
 as follows.
For any locally finite point set $\SS \subset \R^d$, and any
integer $m \in \N$ define
$$
\overline{\xi}(\SS;m) := \sup_{n \in \N} (\essup_{m,n} \{ \xi(\0;
(\SS \cap B(\0;m)) \cup \A) \})
$$
and
$$
\underline{\xi}(\SS;m):= \inf_{n \in \N}( \essinf_{m,n} \{ \xi(
\0; (\SS \cap B(\0;m)) \cup \A) \}),
$$
where $\essup_{m,n}$ (respectively $\essinf_{m,n}$) is essential supremum
(infimum),  with respect to
Lebesgue measure on $\R^{dn}$,
  over sets $\A \subset \R^d \setminus B(\0;m)$
of cardinality $n$.
 Define $\xi_\infty(\SS)$,
   called the {\em limit}
of $\xi$ on $\SS$, by
$$
\xi_\infty(\SS)
 := \limsup_{m \to \infty} \overline{\xi}(\SS;m).
$$
 We shall say the functional $\xi$ {\em stabilizes}
on $\SS$ if
\begin{equation} \label{stab}
\lim_{m \to \infty} \overline{\xi}(\SS;m)
= \lim_{m \to \infty} \underline{\xi}(\SS;m)
= \xi_\infty(\SS).
\end{equation}

For $\tau \in (0,\infty)$, let ${\cal P}_\tau$
 be  a homogeneous Poisson point process of  intensity
$\tau$ on $\R^d $. We are interested particularly in functionals
that stabilize almost surely on  $\Po_\tau$. Note that  with
probability 1,
 $\overline{\xi} (\Po_\tau;m)$
is   nonincreasing in $m$ and
 $\underline{\xi} (\Po_\tau;m)$
is nondecreasing in $m$,
 so they both converge. Stabilization
means they converge to the same limit, almost surely. The present
formulation of stabilization is weaker than that of \cite{PY2}.
Any functional $\xi(x; \X)$ which depends only on the points of
$\X$ within a fixed distance of $x$ is  stabilizing on $\Po_\tau$.


We are interested in functionals  on
 spatial point processes involving non-uniform points, defined as follows.
Let $X_1,X_2,\ldots$ be i.i.d. $d$-dimensional random  variables
 with common  density $f$, which is fixed but  arbitrary.
Define the induced {\em binomial point processes}
\begin{equation} \label{bin}
\X_n:= \X_{n}(f): = \{X_{1},...,X_{n}\}, ~~~ n \in \N.
\end{equation}
Our general limit theory is not for $H_\xi(\X_n)$, but
for $H_{\xi_n}(\X_n)$, where we define
\bea
\xi_n(x;\X) :=  \xi( n^{1/d} x;n^{1/d} \X ).
\label{xin}
\eea

To obtain a LLN for  $H_{\xi_n}(\X_n)$ we use the following approach.
 By coupling $n^{1/d} \X_n$ to a Poisson process of
varying intensity, we show that the local behavior of
$\xi(n^{1/d}X_1; n^{1/d} \X_n)$ is approximated by the local
behavior of the coupled Poisson process.  If the functional $\xi$
stabilizes on homogeneous Poisson point processes, then a
conditioning argument shows distributional  convergence of
$\xi(n^{1/d}X_i; n^{1/d} \X_n)$ for each $X_i \in \X_n$.  Under
appropriate moment conditions on $\xi$, this gives a weak LLN for
$H_{\xi_n}$.  This formalizes the intuitive notion that the
limiting behavior of $H_{\xi_n}$ on finite sets is related to the
behavior of $\xi_\infty$ on the infinite set $\Po_\tau$.

Many of the applications that we consider
are concerned with functionals of graphs
 of the form $G=G(\X)$ defined for each locally finite point set
 $\X \subset \R^d$,
 where either $G(\X)$ or (in the case of the Voronoi graph) its planar dual
has vertex set $\X$.

 We shall say $G$ is {\em translation invariant} if
translation by $y$ is a graph isomorphism from
$G(\X)$ to  $G(y+\X)$ for all $y\in \R^d$ and all locally
finite  point sets $\X$.
We shall say $G$ is {\em scale invariant} if
 scalar multiplication by $a$ induces a graph isomorphism
from $G(\X)$ to
 $G(a\X)$  for all  $\X$ and all $a >0$.

It is useful to have a notion of stabilization for these graphs.
Given $G$, and given a vertex $x \in \X$, let $\Ed(x;G(\X))$ be
the set of edges of $G(\X)$ incident to $x$ (or for the Voronoi
graph, the set of edges whose planar duals are incident to $x$).
Let $\Po_{\tau,0} := \Po_{\tau} \cup \{\0\}$.
 We shall say that $G$
{\em stabilizes} on $\Po_\tau$ if there exists a random but
almost surely finite variable $R$ such that
$$
\Ed(\0;G(\Po_{\tau,0}) ) = \Ed(\0; G(\Po_{\tau,0} \cap B(\0;R) )
\cup \A)
$$
for all finite $\A \subset \R^d \setminus B(\0;R)$.

Stabilization of the graph $G$ says that the local behavior
 of the graph in a bounded region   is unaffected by points
beyond a finite (but random) distance from that region. As we
shall see,  the minimal spanning tree and the
$k$-nearest neighbors,
Voronoi, Delaunay, and sphere of influence graphs are all
stabilizing on $\Po_\tau$, $\tau \in (0, \infty)$.


\subsection{General  LLN results }
\label{gensec}

The following theorem places the objective method
 in a  general context,
 shows that the asymptotic behavior of
 $H_{\xi_n}(\X_n)$ is sensitive to the underlying density $f$, and
  explicitly identifies the asymptotic
constants  in terms of $f$ and the limit functional $\xi_\infty$.
It will be proved in Section \ref{secpfs}.

\begin{theo}\label{generalLLN}
 (General LLN) Suppose $q = 1$ or $q = 2$. Suppose
 $\xi$    is almost surely  stabilizing on $\Po_{\tau}$,  with limit
$\xi_\infty(\Po_{\tau})$, for all $\tau \in (0,\infty)$.
If $\xi$ satisfies the  moments condition
\begin{equation}\label{bpm}
\sup_{ n \in \N} E [ \xi(n^{1/d} X_1; n^{1/d} \X_n )^p  ] <
\infty,
\end{equation}
 for some $p > q$,  then as $n \to \infty$,
 \bea
 \label{mar26}
  n^{-1 } H_{\xi_n}(\X_{n} )
\to \int_{\R^d} E [\xi_\infty (\Po_{f(x)})] f(x) dx \ \ \ {\rm
~in~}L^q.
 \eea
\end{theo}
Obviously,  for there to be any possibility at all for the mean
of   the left side of (\ref{mar26}) to converge to a finite
limit, the  moments condition (\ref{bpm}) must  hold for $p=1$.
In this sense, when $q = 1$,  the moments condition (\ref{bpm})
is close to being the best possible.

A simplification arises in the  case  where
there is a constant $\gamma >0$ such that $\xi$
satisfies the  relation
$$
\xi(a x ; a \X) = a^\gamma \xi(x; \X)
$$
for all positive scalars $a$ and all finite point sets $\X$ and
$x \in \X$. In this case we say $\xi$ is {\em homogeneous of
order $\gamma$}. Homogeneity  of order $\gamma$ implies that
$\xi_n(x;\X)= n^{\gamma/d}\xi(x; \X)$. Moreover, almost sure
stabilization on $\Po_1$ with limit $\xi_\infty(\Po_1)$, together
 with
homogeneity of order $\gamma$, implies
 stabilization on $\Po_\tau$
with  limit $\tau^{-\gamma/d}\xi_\infty(\Po_1)$, for any $\tau
>0$. Therefore the $L^q$ limit in (\ref{mar26}) simplifies to
 \bea
%
 E[ \xi_\infty (\Po_1)] \int_{\R^d}
 f(x)^{(d-\gamma)/d} dx.
 \label{genll}
 \eea

Even simpler is the special case where the function $\xi$ is {\em
scale invariant}, i.e., is homogeneous of order 0. In this case
the expresion (\ref{genll}) simplifies to
 $E[ \xi_\infty (\Po_1)]$, and
the large $n$ behavior of scale invariant functionals is not
sensitive to the density of the underlying point set.

Theorem \ref{generalLLN} admits the following extension to
functionals defined on marked  point sets.   Let
$(\K,\F_\K, P_\K)$ be a probability space.
A {\em marked point set} is a subset of $\R^d \times \K$,
to be denoted $\tilde{\X}$ where $\X$ is an (unmarked) subset
of $\R^d$ and  the tilde indicates that each element $x$ of $\X$ carries
a mark in $\K$ (and with its mark, is denoted $\tilde{x}$).
In this context, a functional $\xi(\tilde{x} ;\tilde{\X})$ is
said to be translation invariant if for all
 $y\in \R^d$, and for any element
$\tilde{x}$ of any marked point set $\tilde{\X}$,
we have $\xi(\tau_y(\tilde{x});\tau_y( \tilde{\X})) =
 \xi(\tilde{x};\tilde{\X})$ where
$\tau_y$ is the translation operator sending any element
$(x,t)\in \R^d \times \K$ to $(y+x,t)$ (i.e., leaving the mark
unchanged).

In the random setting, assume the marks are i.i.d. with
distribution $P_\K$. We are interested in the cases where $\X$ is
the point process $\X_n$ or $\Po_\tau$; in both cases assume the
mark values are independent of the set $\X$. We say that $\xi$
stabilizes on the marked Poisson point process
$\tilde{\Po}_{\tau}$ if (\ref{stab}) holds with $\SS$ replaced by
$\tilde{\Po}_{\tau}$.

If $\xi$ is translation invariant and  almost surely stabilizing
on $\tilde{\Po}_{\tau}$, with limit
$\xi_\infty(\tilde{\Po}_{\tau}), \ \tau \in (0, \infty)$, and if
$ \xi$   satisfies the moments condition (\ref{bpm})
 for some $p > q$, then as $n \to \infty$, we obtain
a version of
 (\ref{mar26}) for marked processes, namely,
 \bea
 \label{markedlimit}
  n^{-1 } H_{\xi_n}(\tilde{\X}_{n} )
\to  \int_{\R^d} E [\xi_\infty (\tilde{\Po}_{f(x)})] f(x) dx \ \
\ {\rm ~in~}L^q. \eea
In applications of (\ref{markedlimit}), it
will be clear that we are considering marked point processes and
for simplicity we will thus suppress mention of the tilde.

Many applications of Theorem  \ref{generalLLN}
are  defined in terms of functionals of graphs
 arising in computational geometry. Suppose the graph
$G:=G(\X)$ is defined for all locally finite $\X$. Given $G$,
functionals such as  total length, number of edges, or
 number of edges  less than some specified length  are of interest.
These and other functionals may be interpreted as a
 total of $\phi$-weighted edge lengths, i.e.,
as a sum \bea L_\phi^G(\X):= \sum_{e \in G(\X)} \phi(|e|),
\label{phiw1} \eea with $\phi: [0,\infty]\to [0,\infty)$ a
specified function. Also of interest are the number of
components, which we denote $K^G(\X)$,  and,
 for any specified finite connected unlabeled graph $\Gamma$,
the number of vertices $x\in \X$  for which $G(\X)$ contains a
subgraph isomorphic to $\Gamma$ with a vertex at $x$, which we
denote $V^G_\Gamma(\X)$.
 Let $\sigma_G $ be the order
 of the component (i.e., the number of vertices in the component)
 containing the origin of
$G(\Po_{1,0})$, and let $E_\Gamma$ be the event that
$G(\Po_{1,0})$ contains a subgraph isomorphic to $\Gamma$ with a
vertex at the origin.

Note that both $K^G(\X)$ and $V_\Gamma^G(\X)$, as well as
$L_\phi^G(\X)$ in the case $\phi \equiv 1$  (the total number of
edges) are scale invariant functionals of $\X$.

The following general result, proved in Section 3, is a
consequence of Theorem \ref{generalLLN}.

\begin{theo}
\label{Gtheo} Suppose $G$ is translation and scale invariant and
stabilizes on $\Po_1$. Then \bea n^{-1} K^G(\X_n) \to
E[\sigma_G^{-1}] \mbox{ in } L^2 \label{Gth1} \eea and for any
finite connected graph  $\Gamma$, \bea n^{-1} V_\Gamma^G (\X_n)
\to P[E_\Gamma] \mbox{ in } L^2. \label{Gth2} \eea Moreover, if
$q=1$ or $q=2$ and
 $\phi: [0,\infty]\to [0,\infty)$ is a specified function
with \bea \sup_{n \in \N} E\left[ \left( \sum_{e \in
\Ed(X_1;G(\X_n)) } \phi(n^{1/d}|e|)
 \right)^p \right] < \infty
\label{Gth2a}
\eea
for some $p>q$,
then
\bea
\label{Gth3}
n^{-1} L_\phi^G( n^{1/d} \X_n) \to
 \frac{1}{2} \int_{\R^d} E
\sum_{e \in \Ed(\0;G(\Po_{1,0}) ) }
\phi \left(  \frac{|e|}{f(x)^{1/d}} \right) f(x) dx
 \ \
{\rm ~in~}L^q.
 \eea
\end{theo}

  The conclusions of Theorems \ref{generalLLN} and \ref{Gtheo}
can be strengthened in more than one way.   If, for example, $q >
2$ is an integer, and if $\xi$ satisfies the moments condition
(\ref{bpm}) for some $p > q$, then a modification of the coupling
arguments given in Section 3 yields convergence in $L^q$.
Moreover, for many functionals, convergence of the means  $E
\left[ n^{-1 } H_{\xi_n}(\X_{n} ) \right]$, as given by Theorems
\ref{generalLLN} and \ref{Gtheo},  implies almost sure and even
complete convergence using concentration inequalities involving
either isoperimetric methods or sums of martingale differences;
we will not pursue this but refer to \cite{JiY, PY2, Y1}
  for details.  We will also not seek to establish rates
of convergence in our general results.

\vskip.3cm

{\em Remarks}.

(a) In the special case where $\phi(x)= x^\alpha$ (i.e.,
power-weighted edges), the left-hand side of (\ref{Gth3}) equals
$n^{(\alpha -d)/d} L_\phi^G(\X_n)$, while the right-hand side of
(\ref{Gth3}) \ simplifies \ to \ $ C_G(d,\alpha)\int_{\R^d}
f(x)^{(d-\alpha)/d}dx$, with $ C_G(d,\alpha):= \frac{1}{2} E
\sum_e \phi(|e|) $, the sum being over all $e \in \Ed(\0;
G(\Po_{1,0}))$. This provides extra information about, for
example, the limiting constant in (\ref{bhh}).

 (b)
The stabilizing hypothesis of Theorem \ref{generalLLN} can be
weakened to one requiring that $\xi$ be almost surely stabilizing
on $\Po_{f(x)}$,  with limit $\xi_\infty(\Po_{f(x)})$, for almost
all $x$ in the support of  the density $f$. These weakened
hypotheses are used in Theorem \ref{thboo} below.  The moments
condition (\ref{bpm}) is not always easy to check, but is
obviously true for any $p$ when the functional $\xi$ is uniformly
bounded. In Section 2.3 we shall verify the moments condition
(\ref{Gth2a}) for various graphs.

(c)
Even without the moments condition (\ref{bpm}),
the stabilization  assumption in Theorem \ref{generalLLN}
is enough to guarantee that the right side of
(\ref{mar26}) is a lower bound for $\liminf E [n^{-1} H_{\xi_n}(\X_n)]$;
this is proved by following the proof of Lemma \ref{lemapr26} below
and applying
 Fatou's Lemma. Likewise, even without the moments condition (\ref{Gth2a}),
a weaker version of (\ref{Gth3}) holds in which the right hand
side is a lower bound for the liminf of the expectation of the
left hand side. With the moments condition (\ref{bpm})
(respectively, (\ref{Gth2a})) the integral on the right hand side
of (\ref{mar26})  (respectively, (\ref{Gth3})) is finite.

(d) Limit laws such as (\ref{Gth1}) and (\ref{Gth2}) for scale invariant
functionals are of interest in multi-dimensional scaling
\cite{Gr}. In this context, given the matrix of interpoint
distances between pairs of points in a data set
$\X$,
one seeks to identify the dimension in which a
data set lives.   To identify the underlying dimension, it is
useful to study scale invariant functionals of the interpoint
distances, since these are precisely the functionals whose
asymptotics are sensitive only to the dimension of the support of
the distribution of the data points and not on the underlying
density. This approach to dimension identification may be
relatively inexpensive from a computational point of view
\cite{BQY}.

(e) Without further conditions on $\xi$ and $H_{\xi_n}$, we are
unable to obtain asymptotics for $H_{\xi_n}$ over point sets
consisting of random $d$-vectors having a law with a singular
component.

(f) Suppose the conditions of Theorem \ref{generalLLN} are
satisfied with $q = 2$. Since $\xi^2$ is almost surely
stabilizing on $\Po_{\tau}$,  with limit
$\xi^2_\infty(\Po_{\tau})$, and since $\xi^2$ satisfies the
moments condition (\ref{bpm}) for some $p > 1$, it follows that
the sample variance of $\{\xi_n(X_i;\X_n), 1 \leq i \leq n\}$,
 namely the quantity
$$
\frac { \sum_{i=1}^n (\xi_n(X_i; \X_n)) - n^{-1} \sum_{i=1}^n
\xi_n(X_i; \X_n))^2 } {n}
$$
converges in $L^1$ to $\text{Var}(\xi_{\infty} ( \Po_{f(X_1)}) ).$

(g) The limit (\ref{Gth2}) says, loosely speaking, that the number
of vertices of $G(\X_n)$ satisfying any property
determined by the local graph structure within a bounded
graph distance exhibits LLN behavior. For example,
 it    yields a LLN for the proportionate
number of vertices of any fixed degree, among other things.

(h) A version of Theorem \ref{Gtheo} also holds for directed
graphs. In this context the limit (\ref{Gth3}) holds without the
factor of $1/2$  and with ${\cal E}(x;G(\X))$ defined to be the set
of edges going into $x$.

\subsection{Applications}
\label{secapp}

The applications of Theorems \ref{generalLLN} and \ref{Gtheo}
range from the treatment of functionals in computational geometry
to the statistics of Boolean models as well as to an analysis of
packing processes.  The following discussion is not meant to be an
exhaustive treatment of applications, but is merely meant to
indicate the variety of uses of the main theorems.


\vskip.5cm

(a) {\em Minimal spanning tree}. Given a locally finite set $\X
\subset \R^d, \ d \geq 2$, let $\MST(\X)$ be the graph with
vertex set $x$ obtained by including each edge $e:=\{x,y\}$ for
which  there is no  path in $\X$ from $x$ to $y$ consisting of
links which are all shorter than $e$. If $\X$ is finite with
distinct inter-point distances, then $\MST(\X)$  is the
  {\em minimal spanning tree} on $\X$, i.e. the
connected graph with vertex set $\X$ of minimal total edge length;
see Alexander \cite{Ale}. If $\X$ is infinite, then $\MST(\X)$ is
the so-called {\em minimal spanning forest} on $\X$ \cite{Ale}.
Clearly $\MST(\X)$ is translation and scale invariant.

The following theorem  considerably expands upon Theorem
\ref{basicMST}. We say that $\phi$ has {\em polynomial growth of
order} $a < \infty$ if $\phi(x) \leq C(1 + x^a)$.


\begin{theo}\label{MST}
Let $G(\X) = \MST(\X)$. Then (\ref{Gth2}) holds for any finite connected
$\Gamma$. Also,
(\ref{Gth3}) holds with $q=2$ if
 either (i) $\phi$ is bounded, or
 (ii) $\phi$ has polynomial growth of order $a < \infty$,
 the support of $f$ is a convex polyhedron and  $f$ is
bounded away from infinity and zero on its support, and
(\ref{Gth3}) holds with $q=1$ if (iii) $\phi$ has polynomial
growth of order $a < d$ and $\int_{\R^d} \vert x \vert^r f(x) dx
< \infty$, for some $r > \max \{ad/(d-a), d/(d-a) \}$, and
$\int_{\R^d} f(x)^{(d-a)/d}dx < \infty$.
\end{theo}

Special cases of Theorem \ref{MST} include an $L^2$ version of
 Theorem \ref{basicMST}
(take $\phi$ to be a power function)
and the $L^2$ LLN for the empirical distribution function of
rescaled edge
 lengths in the MST, first obtained by Penrose \cite{P1}
(take  $\phi(x):= 1_{[0,t]}(x)$). For uniformly distributed
points, Bezuidenhout et al. \cite{BGL} obtained an almost sure
version of the latter result. The first part (\ref{Gth2}) of the
conclusion in Theorem \ref{MST} yields an $L^2$ version of the
results of \cite{SSE} and Theorem 5.2.2 of \cite{Stbk}, concerning
the number  of vertices of fixed degree.

\vskip.3cm

{\em Proof of Theorem \ref{MST}.} The set of edges
of $\MST(\Po_{1,0})$ incident to the origin is unaffected by any
additions or deletions of points outside a ball of random but
almost surely finite radius, i.e.   $G(\X)= \MST(\X)$ is
stabilizing on $\Po_1$.
 This follows from the definition of $\MST(\X)$ and
the simultaneous uniqueness of percolation clusters (Alexander
\cite{Ale2}), and the almost sure existence of a ``blocking set''
of Poisson points lying in some annulus surrounding the origin and
precluding the possibility of any edge  from  the origin to the
exterior of that annulus. For more details see Lee \cite{Lee}.

Theorem \ref{Gtheo} then  gives us (\ref{Gth2}).
To prove (\ref{Gth3}) it remains only to verify the moments
condition (\ref{Gth2a}).
 Note first that points in the MST
 have a degree which is
uniformly bounded by some finite constant $C(d)$ (Lemma 4 of
\cite{AS}).
 Hence under condition (i), i.e. if $\phi$ is bounded,
then $
 \sum_{e \in \Ed(x;G(\X)) } \phi(|e|)
$ is also bounded for any point set $\X$ and $x \in \X$, and thus
(\ref{Gth2a}) is clearly satisfied, showing that (\ref{Gth3})
holds in this case.

The uniform bound on vertex degrees implies that for any $p>1$,
there is a second constant $C(p,d)$  such that for any $\X$ and
any $x \in \X$, \bea
 \left( \sum_{e \in \Ed(x;G(\X)) } \phi(n^{1/d}|e|) \right)^p
 \leq C(p,d) \sum_{e \in \Ed(x; G(\X)) }
\phi^p(n^{1/d}\vert e \vert). \label{jun} \eea

 Suppose condition
(ii) holds, i.e. $\phi$ has polynomial  growth of order $a<
\infty$ and the support of $f$ is a convex polyhedron with  the
restriction of $f$ to its support bounded away from zero and
infinity. Then  by the uniform bound on vertex degrees, the right
hand side of (\ref{jun}) is bounded by a constant plus a constant
multiple of \bea
 n^{ap/d} \max_{e \in \Ed(x;G(\X))} |e|^{ap}.
\label{may7} \eea  Lemma 2.1 of \cite{Y2} shows that
$$
E[\max \{ |e|^{ap}: e\in \Ed(X_1;G(\X_n))\} ] = O(n^{-ap/d})
$$
and therefore the left side of (\ref{Gth2a})
  is uniformly bounded by a constant.
Therefore, under condition (ii), the conclusion (\ref{Gth3}) (with
$q=2$) holds.

Finally consider case (iii). Put $\xi_n(x;\X):= \sum_{e \in
\Ed(x,G(\X)) } \phi(n^{1/d}|e|)  $.
 If $\phi$ has polynomial growth of order $a$ then
 by exchangeability and (\ref{jun}),
\bea
 E [(\xi_{n}(X_1;  \X_n))^p ]
 = n^{-1}
E\sum_{i=1}^n
(\xi_{n}(X_i;  \X_n))^p \nonumber \\
 \leq C(p,d) n^{-1}E \sum_{i=1}^n \sum_{e \in  \Ed(X_i;G(\X_n))}
\phi^p( n^{1/d} \vert e \vert)
\nonumber \\
 =2 C(p,d) n^{-1}E \sum_{e \in G(\X_n)} \phi^p( n^{1/d} \vert e \vert)
\nonumber \\
\leq C_1(p,d) + C_2(p,d) n^{(ap - d)/d}  E \sum_{e\in G(\X_n)}
\vert e \vert^{ap}. \label{jun20} \eea

We now bound the right hand side of (\ref{jun20}). Let $L^{ap}
(\X) := \sum_{e \in G(\X)} \vert e \vert^{ap} $ and for all $k
\geq 1$, let $A_k$ denote the annular shell centered around the
origin of $\R^d$ with inner radius $2^k$ and outer radius
$2^{k+1}$ (Let $A_0$ be the ball centered at the origin with
radius $2$). Note that, as in (7.21) of \cite{Y1},
\begin{equation}\label{7.21}
 L^{ap}(\X_n) \leq
\sum_{0 \leq k \leq s(n)} L^{ap}(\X_n \cap A_k) +
    C(p) \max_{1
\leq i \leq n} \vert X_i \vert^{ap},
\end{equation}
 where $s(n)$ is the largest
$k$ such that $A_k \cap \X_n$ is not empty. We need to show,
after taking expectations and dividing by $n^{(d - ap)/d}$, that
the two terms on the right are bounded uniformly in $n$.   By
Jensen's inequality and the growth bounds
 $L^{ap}(\X) \leq C
(\text{diam} \X)^{ap} (\text{card} (\X))^{(d-ap)/d}$
(see Lemma 3.3 of \cite{Y1}),
 the first
term is upper bounded by
$$
C \sum_{k \geq 1} 2^{kap} ( P [X_1 \in A_k])^{(d-ap)/d}
$$
which is finite by the integrability hypothesis (see p. 85 of
\cite{P1}).  The second term is  bounded by
\bea
C(p) \int_0^{\infty}  P\left[ \max_{1 \leq i \leq n} \vert X_1 \vert^{ap}
 \geq t
n^{(d-ap)/d} \right]  dt
\\
\leq  n C(p)
\int_0^{\infty}  P\left[  \vert X_1 \vert^{apd/(d-ap)}
 \geq t^{d/(d-ap)}
n \right]  dt. \eea
By  Markov's inequality
together with the moment condition $\int_{\R^d} \vert x \vert^r
f(x) dx < \infty$, for some $r > ad/(d-a)$, this last integral is
finite.
 Thus under condition (iii),
we can choose $p > 1$ so that (\ref{jun20}) is uniformly bounded,
and hence the moments condition (\ref{Gth2a}) is satisfied,
showing the validity of the conclusion under condition  (iii).
 $\qed$ \\

\vskip.3cm

 (b) {\em $k$-nearest neighbors graphs}.
Let $k$ be  a positive integer.
Given a locally finite  point set $\X \subset \R^d$,
  the $k$-nearest neighbors (undirected) graph
  on $\X$, denoted $\NG(\X)$, is the graph with vertex set
$\X$ obtained by including $\{x,y\}$ as an edge whenever $y$ is
one of the $k$ nearest neighbors of $x$ and/or $x$ is one of the
$k$ nearest neighbors of $y$. The $k$-nearest neighbors
(directed) graph on $\X$, denoted $NG'(\X)$, is the graph with
vertex set $\X$ obtained by placing an edge between each point
and its $k$ nearest neighbors.
 If the $k$-th nearest neighbor of $x$ is not
well-defined (i.e., if there is a ``tie" in the ordering of
interpoint distances involving $x$), use the lexicographic
ordering as a ``tie-breaker" to determine the $k$ nearest
neighbors. Such a tie has zero probability  for  the random point
sets under consideration here.

It is clear that the $k$-nearest neighbors graphs are translation
and scale invariant. The following generalizes and extends the
asymptotics for the sum of power-weighted edge lengths by
McGivney \cite{M},
 Yukich (Theorem 8.3 of \cite{Y1}), and Jimenez and Yukich
 \cite{JiY},
 who limit attention to continuous densities and increasing
 $\phi$.
It also extends Theorem 2 of Eppstein, Paterson and Yao
\cite{EPY}, who prove convergence of the mean   number of
components for uniform samples.

\begin{theo}
\label{knn} Let $G(\X)$ denote either $ \NG(\X)$ or $\NG'(\X)$.
 Then (\ref{Gth1}) and (\ref{Gth2}) hold.
Moreover, (\ref{Gth3}) holds, with $q =2$,
  if either
(i) $\phi$ is bounded, or
 (ii) $\phi$ has polynomial growth of
 order $a<\infty$, the support of $f$ is a convex polyhedron and  $f$ is
bounded away from infinity and zero.  Finally the directed graph
version of
 (\ref{Gth3}) holds for $\NG'(\X)$ with $q =1$
  if  $\phi$ has polynomial growth
of order $a < d$, $\int_{\R^d} f(x)^{(d-a)/d}dx< \infty$, and
$\int_{\R^d} \vert x \vert^r f(x) dx < \infty$ for some $r >
d/(d-a)$.
\end{theo}

{\em Proof.}
We apply Theorem \ref{Gtheo}.
As shown in  Lemma 6.1 of \cite{PY1} (even though the definition
of stabilization there is slightly different),
the set of edges incident to  the
 origin in $\NG(\P_{1,0})$ is unaffected by
the addition or removal of points  outside a ball of random but
almost surely finite radius, i.e. the graph $G(\X)= \NG(\X)$ is
 stabilizing on $\Po_1$. Similar arguments show that
 $\NG'(\X)$ is
 stabilizing on $\Po_1$.
 Therefore Theorem \ref{Gtheo} yields
(\ref{Gth1}) and (\ref{Gth2}) for $\NG(\X)$ and $\NG'(\X)$.

To prove (\ref{Gth3}) it remains only to verify the moments
condition (\ref{Gth2a}).  This condition is checked, under any of
the conditions (i), (ii), or (iii), in very much the same manner
as for the MST of Theorem \ref{MST}.  The existence of a uniform
bound on the degree of vertices in the nearest neighbors graph is
Lemma 8.4 of \cite{Y1}. Under condition (iii), we note that since
$\NG'$ is subadditive without any error term, the last term in
(\ref{7.21}) is not needed, eliminating the need for the
condition $r > ad/(d-a)$.
$ \qed$

\vskip.3cm

 Henze \cite{He} considers the fraction of points in
$\X_n$
 which are the nearest neighbors of
exactly $j$ other points; he also considers the fraction of points
in $\X_n$ that are the $l$th nearest neighbors to their own $k$th
nearest neighbor.
  He shows that the limits of these fractions
converge to explicit but rather complicated limiting constants.
A directed-graph version of (\ref{Gth2}) in Theorem \ref{Gtheo}
 yields Henze's results  and shows that one
can interpret his limiting constants in terms of a functional
evaluated at a point in the origin of the homogeneous Poisson
point process.

 \vskip.3cm

\vskip.3cm

(c) {\em Voronoi and Delaunay graphs}. Given a locally finite
 set $\X \subset \R^d$, and given $x \in \X$,  the locus of points
closer to $x$ than to any other point $ \in \X$ is called the {\em
Voronoi  cell} centered at $x$. The graph on vertex set $\X$ in
which each pair of adjacent cell centers is connected by an edge
is called the {\em Delaunay graph} on $\X$; if $d=2$, then the
planar dual graph consisting of all boundaries of Voronoi cells
is called {\em Voronoi graph} generated by $\X$.
Edges of the Voronoi graph can be finite or infinite.
Let $\Del(\X)$ (respectively $\Vor(\X)$)
 denote the collection of edges in the Delaunay graph
 (respectively,
the Voronoi graph)  on $\X$.
 The Voronoi and Delaunay graphs
are clearly scale and translation invariant.

\begin{theo}\label{vor} Let $d=2$ and let $G(\X)= \Vor(\X).$
 Then (\ref{Gth2}) holds. Also, if $\phi$
has polynomial growth with $\phi(\infty)=0$,
and if the support of $f$ is a convex
polygonal region and if $f$ is bounded away from infinity and
zero on its support, then (\ref{Gth3}) holds with $q=2$.
\end{theo}

This  result adds to McGivney and Yukich \cite{MY} and Jimenez
and Yukich \cite{JiY}, who require  continuous $f$ and functions
$\phi$ which are either power functions or increasing. If $\phi$
is the identity and the support of $f$ is the unit square, then
the right hand side of (\ref{Gth3}) simply reduces to $2
\int_{[0,1]^2} f(x)^{1/2} dx$ (see Theorem 1.2 of \cite{MY}).

{\em Proof.} Once again we  apply Theorem \ref{Gtheo}. As in
\cite{PY1}, we can verify $G(\X)$ is stabilizing, i.e., the
Voronoi cell centered at the origin for $\Po_{1,0}$  is
unaffected by changes beyond a random but almost surely finite
distance from of the origin. Let $\Ed(X;\Vor(\X))$ denote the
edges of the cell around $X_1$ in $\Vor(\X_n)$.
 Concerning the moments condition
(\ref{Gth2a}), a modification of Lemma 8.1 of \cite{PY1} shows
that under the prescribed conditions on $f$ and its support,
\begin{equation}\label{bpm1}
\sup_{ n \in \N} E \left[  \left( \sum_{e \in
\Ed(X_1;\Vor(\X_n))}  (n^{1/d}|e|)^a \right)^p \right] < \infty
\end{equation}
for all $a > 0$ and $p > 2$. A modification of the proof of Lemma
2.6 of \cite{MY} (pp. 286-87) shows that the cardinality of
$\Ed(X_1;\Vor(\X_n))$ has a finite $p$th  moment for all $p > 2$.
Combining this with (\ref{bpm1}) thus shows that for all $\phi$
with polynomial growth
\begin{equation}
\sup_{ n \in \N} E \left[  \left( \sum_{e \in
\Ed(X_1;\Vor(\X_n))} \phi(n^{1/d}|e|) \right)^p \right] < \infty
\end{equation}
for all  $p > 2$. Thus (\ref{Gth2a}) holds and the conclusion
(\ref{Gth3}) follows.
$\qed$

A similar result to Theorem \ref{vor} holds
if $G(\X)$ is taken to be
 the Delaunay graph and $\Ed(x;\Del(\X))$ denotes the edges in
 $\Del(\X)$ incident to $x$.
 In this case  the conclusion (\ref{Gth2})
yields a LLN for the number  of vertices of the Delaunay graph
of fixed degree $m$,
  for any $m = 3, 4, 5,...$.
Since this quantity is the same as
the total number of cells in the Voronoi graph on $\X_n$ which are $m$-gons,
 this adds to results of Hayen and
 Quine \cite{HQ} who determine the proportion of triangles in the
 Voronoi graph on ${\cal P}_1$.
\vskip.3cm

Finally, for each $t>0$ consider the case
where $\xi(x; \X)$ equals $1$ or $0$ according to whether
 the area of the Voronoi cell around $x$ is bounded by $t$ or not.
This is one case where it is  natural  to use Theorem \ref{generalLLN}
rather than Theorem \ref{Gtheo}, and that result yields a LLN for
 $H_{\xi_n}(\X_n):= \sum_{i=1}^n \xi(n^{1/d}X_i; n^{1/d} \X_n)$,
i.e., a LLN for the empirical distribution function of
 the rescaled areas of the Voronoi diagram on $\X_n$.

 \vskip.3cm

 (d)  {\em Sphere of influence graph}.
Given a locally finite set  $\X \subset \R^d$,
 the sphere of influence graph $\SIG(\X)$ is
a graph with vertex set $\X$, constructed
as follows:  for each $x \in \X$ let $B_x$ be a ball around $x$ with
radius equal to $\min_{y \in \X \setminus \{x\}}
 \{\vert y - x \vert \}.$  Then $B_x$ is
called the sphere of influence of $x$.  Draw an edge between
$x$ and $y$ iff the balls $B_x$ and $B_y$ overlap.
The  collection of such edges is the sphere of influence graph (SIG) on $\X$.
It is clearly translation and scale invariant.

The following LLN is apparently new, even for the identity
function $\phi(x) = x$.
 In the case
$\phi(x)\equiv 1$ it extends  a result of
F\"uredi \cite{Fu} on  the mean number of edges of the SIG
on uniform point sets
 (F\"uredi identified the limiting constant in this case).

\begin{theo}\label{SIG} Let $G(\X)= \SIG(\X)$.
Then (\ref{Gth1}) and (\ref{Gth2}) both hold. If $\phi$ has
polynomial growth and if the support of $f$ is a convex
polyhedron and if $f$ is bounded away from infinity and zero on
its support, then (\ref{Gth3}) holds with $q = 2$.
\end{theo}
{\em Proof.} We apply Theorem \ref{Gtheo}. As in  \cite{PY1}, we
can check that the edges incident to the origin in
$\SIG(\Po_{1,0})$ are  unaffected by  changes beyond a
 random but
almost surely finite distance of the origin, i.e. the graph
$G(\X)=\SIG(\X)$ is stabilizing. Concerning the moments condition
(\ref{Gth2a}), the arguments of \cite{PY1} (Theorem 7.2) show that
for any $a > 0, p>1,$ we have under the prescribed conditions on
$f$,
\begin{equation}\label{bpm2}
\sup_{ n \in \N} E \left[ \left( \sum_{e \in \Ed(X_1; \SIG(\X_n))}
(n^{1/d} |e|)^a \right)^p \right] < \infty.
\end{equation}

Moreover, since the third moment of the degree of  vertices in the
SIG on $\X_n$ is uniformly bounded over all vertices (see e.g. pp.
142-43 of \cite{HJY}) we obtain for any $\phi$ with polynomial
growth
\begin{equation}\label{bpm2a}
\sup_{ n \in \N} E \left[ \left( \sum_{e \in \Ed(X_1; \SIG(\X_n))}
\phi(n^{1/d} |e|) \right)^p \right] < \infty.
\end{equation}
Thus  $(\ref{Gth2a})$ holds.  $\qed$

\vskip.3cm

(e)  {\em Proximity graphs.} Devroye \cite{Dev} defines a
proximity graph on $\X$ to be one in which each $\{x,y\}$ is
included as an edge if a specified set $S(x,y)$ is empty. If
$S(x,y)$ is the ball with opposite poles at $x,y$, then the
associated proximity graph is the {\em Gabriel graph}. If
$S(x,y)$ is the intersection of $B(x;|y-x|)$ and $B(y;|y-x|)$
then it is the {\em relative neighborhood graph}. For a survey of
applications of proximity graphs such as these, and also of the
sphere of influence graph, see Jaromczyk and Toussaint \cite{JT}.

Both the Gabriel graph and the
 relative neighborhood graph are translation and scale invariant,
 and stabilize on
$\Po_1$, and also satisfy (\ref{Gth2a}), subject to conditions on
$\phi$ and $f$ similar to those already given for the sphere of
influence graph; see remarks in \cite{PY1}, Section 9. Therefore
Theorem \ref{Gtheo} yields information about these graphs, adding
to results in Devroye \cite{Dev} on the expected number of edges.
Numerical estimates for values of mean edge length and mean
degree of various proximity graphs over homogeneous Poisson point
sets are given by Smith (chapter III.C of \cite{Sm}).

\vskip.3cm

(f) {\em Boolean models.} As already indicated, Theorem
\ref{generalLLN} extends to marked
 processes via the limit (\ref{markedlimit}). An example
of application is to Boolean models, whose
  importance in stochastic geometry and spatial statistics can
 be seen from e.g.  Hall  \cite{Hall} and  Molchanov \cite{Mo}.
Let $\mu_S$ be a  {\em shape distribution}, that is, a probability
distribution  on the space ${\cal S}$ of all compact sets in
$\R^d$. For measure-theoretic details see Matheron \cite{Math},
page 27. Assume that $\mu_S$ is concentrated on connected sets
contained in $B(\0;K)$ for some fixed finite $K$ (i.e., uniformly
bounded connected sets). On a suitable probability space let
$(S_i, \ i \geq 1)$ be a family of random closed sets each with
distribution $\mu_S$,
 independent of each other and of $(X_1,X_2,\ldots)$
(as usual, $X_1,X_2,\ldots$ are i.i.d. $d$-vectors with
common density $f$).
Let $\Xi_n := \cup_{i = 1}^n(X_i + n^{-1/d}S_i)$.
We refer to  $X_i  + n^{-1/d}S_i  $ as a
 {\em random shape centered at $X_i$}.
The random set $\Xi_n$ is a scale-changed Boolean model in the
sense of Hall \cite{Hall}, pages 141 and 233.

 A connected
component of $\Xi_n$ is often called a {\em clump}.
 A clump of order $k$ is one which
comprises precisely $k$ random shapes. Let $U_k(\Xi_n)$ be the
number of clumps  of order $k$ and let $U(\Xi_n) := \sum_k
U_k(\Xi_n)$, the total number of clumps. Let $V(\Xi_n)$ denote the
total volume of the set $\Xi_n$. In the case $d=2$, we consider
 the {\em total curvature functional}  of the set $\Xi_n$
(the product of $2\pi$ and the Euler characteristic of $\Xi_n$;
see Hall (\cite{Hall}, Ch. 4.3)), which we shall denote
$W(\Xi_n)$.

Another statistic of interest is the {\em off-line packing
functional} for the collection $X_i + n^{-1/d}S_i, \ 1 \leq i
\leq \ n$. This functional, denoted by $M(\Xi_{ n})$, is the
maximal number of non-intersecting random shapes in the
 collection $X_i + n^{-1/d}S_i, \ 1 \leq i \leq  n$.
Additional statistics associated with Boolean models are found in
Molchanov \cite{Mo}, for example.

Let $\Xi_{\infty,\lambda}$  be the infinite Boolean model
 $\cup_{X \in \Po_{\lambda}}
(X+S_X)$, where $\Po_{\lambda}$ is a homogeneous Poisson process
of intensity $\lambda$ on $\R^d$, each point carrying an
independent ${\cal S}$-valued mark $S_X$ with distribution
$\mu_S$.
 Let $\lambda_c := \lambda_c(d,\mu_{S})$ denote the continuum
 percolation threshold,
i.e., let $\lambda_c $  be the supremum of the set of values of
$\lambda $ such that $\Xi_{\infty,\lambda}$ almost surely has no
infinite connected component
 (note $\lambda_c= \infty$ if $d=1$).

\begin{theo}
\label{thboo} (i) There exist constants $u_{k,\infty}$ ($k\in
\N$) and $u_\infty$ (dependent on $\mu_S$), such that
$n^{-1}U_k(\Xi_n) \to u_{k,\infty}$ in $L^2$ (for each $k \in
\N$) and such that $n^{-1}U(\Xi_n) \to u_{\infty}$ in $L^2$.

(ii) There exists a constant $v_{\infty}$  (dependent on $\mu_S$)
such that $V(\Xi_n) \to v_{\infty}$ in $L^2$.

(iii)
If  $d=2$ and the measure
$\mu_S$ is concentrated on convex sets,
then there exists a constant $w_{\infty}$
(dependent on $\mu_S$) such that
 $W(\Xi_n)$ converges in   $L^2 $ to $w_\infty$.

 (iv)
  If $\sup_{x \in \R^d} f(x) < \lambda_c$, then there is a constant
$m_{\infty}$ such that $n^{-1} M(\Xi_{ n})$ converges in $L^2$
to  $ m_{\infty}$.

\end{theo}

Part (i) of this result adds
 to existing results in the literature such as Hall's
result (\cite{Hall},  Theorem 4.7) on the number of clumps of
order 1, especially
since we do not restrict attention to a  uniform underlying
density $f$ for the points $X_i$.
As usual,  the value of the
limiting constant $u_\infty$ is of the form $\int
E[\xi_\infty(\Po_{f(x)})] f(x) dx,$ where  $\xi(x;\X)$ is
described in the proof below.

Part (ii) shows that the volume
functional satisfies a weak LLN over non-uniform point sets,
adding to results of Hall (\cite{Hall}, Ch. 3.4) involving the
vacancy functional of $\Xi_n$. Part (iii)
 is a
weak LLN for the total curvature functional over non-uniform
samples, and also adds to results of Hall (\cite{Hall}, Ch. 4.3)).

  In part (iii), the
constant $w_\infty$ is again of the form
 $\int E[\xi_\infty(\Po_{f(x)})] f(x) dx,$
and if $\mu_S$ is isotropic,  there  exist analytic formulae for
 $E[\xi_\infty(\Po_{\lambda})] $; see (2.14) of \cite{Mo},
or (4.28) of \cite{Hall}.

Concerning part (iv),  the off-line packing functional
 can be shown to be subadditive, and methods based
on this fact show that a weak LLN also holds, at least for
the uniform distribution on a cube of volume $1/\lambda$, even
in the supercritical case $\lambda \geq  \lambda_c$.

{\em Sketch of proof of Theorem \ref{thboo}.} Define the rescaled
Boolean model $\Xi'_n$ by
$$
 \Xi'_n:= n^{1/d} \Xi_n =\cup_{i=1}^n (n^{1/d}X_i + S_i) .
$$
Let the mark space be ${\cal S}$ with mark distribution $\mu_S$.

(i) Given a marked point set $\X \subset \R^d$ with marks $S_x, x
\in \X$, let
 $\xi(x;\X)$ be the reciprocal of
the order of the clump of
 $\cup_{y \in \X} (y + S_y)$
 containing $x+ S_x$.
Then $H_\xi(\X)$ is the total number of clumps of
 $\cup_{y \in \X} (y + S_y)$; hence
$$
H_{\xi_n}(\X_n) = U(\Xi'_n) = U(\Xi_n).
$$
Stabilization of $\xi$ follows from the fact that $\mu_S$ is
concentrated on uniformly bounded sets. The moments condition
(\ref{bpm}) follows from the uniform bound $\xi(x;\X)\leq 1$.
Therefore the LLN for $U (\Xi_n)$ follows from
(\ref{markedlimit}). The LLN for $U_k(\Xi_n)$ is proved similarly.

(ii) This time let $\xi(x;\X)$ be the volume of the intersection
of $\cup_{y \in \X}S_y$
 with the Voronoi
cell around $x$ for $\X$. Then
$$
H_{\xi_n}(\X_n) = U(\Xi'_n) = n U(\Xi_n).
$$
Since $\mu_S$ is concentrated on sets contained in $B(\0;K)$, for
any $x \in \X$, the intersection of
 $\cup_{y \in \X}S_y$ with
the Voronoi cell around $x$ is contained in $B(x;2K)$, since any
point lying outside $B(x;2K)$ but in
 $\cup_{y \in \X}S_y$ must be closer to some point $y\in\X
\setminus \{x\}$ than it is to $x$. Both the stabilization and
the moments conditions in the marked point process version of
Theorem \ref{generalLLN} hold as a consequence of this, and
(\ref{markedlimit}) yields the LLN for $V(\Xi_n)$.

(iii) This time let $\xi(x;\X)$ be the  contribution of the
random set $S_x$ to the total curvature of the
union $\cup_{y \in \X}\Xi_y$. This gives us $H_\xi(\X_n) =
W(\Xi_n)$. By convexity, $\xi(x;\X)$
 is uniformly bounded, which gives us the moments condition
(\ref{bpm}).

By the
 assumption of
 uniform boundedness of the random sets, if
${\cal S}_\0$ is an independent random shape inserted at the
origin,
changes outside $B(\0;2K)$ do not affect the contribution of
${\cal S}_\0$ to the total curvature.
 Thus, for  all $\tau$, $\xi$ is almost surely
stabilizing on $\Po_{\tau}$.
 This enables us to deduce the
stabilization condition in the marked point process version of
Theorem \ref{generalLLN}, and,  we
therefore  deduce the LLN behavior.

(iv)  Let $\xi(x;\X)$ be either $1$ or $0$, depending upon
whether or not $x + S_x$ is included in the maximal subset of
non-intersecting shapes.
 If there are several such
maximal subsets, choose one in an arbitrary
deterministic but translation-invariant manner.
Then $H_{\xi}(\X)$ is the maximal number of disjoint shapes and
$H_{\xi_n} (\X_{ n} ) = M(\Xi'_{ n})$.  Since $\xi \leq 1$, this
gives  the moments condition (\ref{bpm}).

 If $f < \lambda_c$, then there is
almost surely no infinite cluster in $\Po_{f(x)}$, for we are in
the subcritical phase of continuum percolation.  Since there is no
infinite cluster, inserting a random shape at the origin thus
almost surely  changes the order of only finite clusters and thus
the packing functional $M$ stabilizes.
 \qed

\vskip.3cm

(g) {\em Packing processes}. Consider the following prototypical
random sequential packing model. Let
$B_{n,1},B_{n,2}...,B_{n,n}$   be  a sequence
 of $d$-dimensional balls of volume $n^{-1}$  whose centers
are independent random $d$-vectors with
 common
 probability density
function $f: \R^d \to [0, \infty)$. Let the first ball $B_{n,1}$
be {\em packed}, and recursively for $i=2,3, \ldots,n$, let
 the $i$-th ball
$B_{n,i}$ be  packed iff $B_{n,i}$ does not overlap any ball in
$B_{n,1},...,B_{n,i-1}$ which has already been packed. If not
packed, the $i$-th ball is discarded. Let $ N_{f}(n)$ be the
number of packed balls, out of the first $n$ to arrive.
This is sometimes called `on-line packing', in contrast
with the off-line scheme described earlier.

We may use our general result (\ref{markedlimit}) for marked processes
to obtain LLN for random sequential
packing, as follows.

For any finite point set  $\X \subset \R^d$,
assume the points have marks which are independent and
uniformly distributed over $[0,1]$.
Assume unit volume balls centered at the points of
$\X$  arrive sequentially in an order determined by the associated marks,
and assume as before that each ball is packed or discarded according to
whether or not it overlaps a previously packed ball.  Let
$\xi(x; \X)$ be
either $1$ or $0$ depending on whether the ball
centered at $x$ is packed or discarded. Then with
the binomial point process
$\X_n$ defined at (\ref{bin}), and $\xi_n$ defined at
(\ref{xin}), it can be seen that
$H_{\xi_n}(\X_n)$ has the same distribution as  $N_f(n)$.
 Following \cite{PY2}, we can show that $\xi$ is almost
surely stabilizing on $\Po_{\tau}, \ \tau \in (0,\infty)$, with
limit $\xi_{\infty}$. Since $\xi$ is bounded it satisfies the
moments condition (\ref{bpm}) and therefore using the limit
(\ref{markedlimit}) we get the following    LLN
 for $N_{f}(n)$.

\begin{theo}\label{markedLLN} Let $f: \R^d \to [0,\infty) $ be an arbitrary
density.  As $n \to \infty$,
 \bea
\label{jun11}
  n^{-1 } N_{f}(n)
\to \int_{\R^d} E [\xi_\infty (\Po_{f(x)})] f(x) dx \ \ \ {\rm
~in~}L^2.  \eea
\end{theo}
This result represents a finite input version of Theorem 5.1 of
Penrose \cite{P2}, but with a more general class of densities $f$.

Theorem \ref{markedLLN} extends to more general versions of the
prototypical packing model. For example, by following the general
stabilization analysis of \cite{PY2},  we can develop asymptotics
in the finite input setting for the number of packed balls in the
following general models: (i) models with balls replaced by
particles of random size/shape/charge, (ii) time dependent,
dynamic models, (iii) cooperative sequential adsorption models,
and (iv) ballistic deposition models. In each case, we obtain a
LLN for the number of packed balls, from among the first $n$ to
arrive,  when the distribution of the balls has a density $f: \R^d
\to [0,\infty)$.
 See \cite{PY2} for a
discussion of these models and for laws of large numbers in the
special case  where the particless arrive {\em uniformly} at random
over $\R^d$. \\

\vskip.3cm

 (h) {\em Combinatorial optimization}.
Both the shortest traveling salesman tour and the minimal matching are
translation and scale invariant graphs on finite point sets.
However, it is not known whether their definition can be
extended to infinite  sets in a manner that makes them stabilizing,
 and therefore
we are at present unable to apply Theorem \ref{generalLLN} or
\ref{Gtheo} in these cases.

\section{Proofs}
\label{secpfs}
\allco
 The proof of Theorem \ref{generalLLN} centers around suitably coupling 
a version of the
binomial process $\X_n$ to a Cox process, that is,  a Poisson
process whose intensity measure is itself random. We do this as
follows. On a suitable probability space suppose we have,
independently, a $d$-dimensional variable $X$ with density $f$,
and a homogeneous Poisson processes $\Po$ of rate $1$
on $\R^d \times [0,\infty)$.

Define coupled point processes
 $\Po(n)$, $\X'_{n-1}$, and $\H_n$ and a random variable
$\zeta_n$, all in terms of $\Po$ and $X$,
 as follows.
Let $\Po(n)$ be the image of the  restriction of $\Po$ to the set
$$
\{(x,t) \in \R^d \times [0\ ,\infty): t \leq n f(x) \},
$$
 under the projection $(x,t) \mapsto x$. Then $\Po(n)$ is a  Poisson
process in $\R^d$
   with intensity function $nf(\cdot)$, consisting
of  $N(n)$ random points
with common density $f$.
 Discard  $(N(n)-(n-1))^+$ of these  points, chosen at random,
and add $(n-1-N(n))^+$ extra independent points with common density $f$.
The resulting set of points is denoted $\X'_{n-1}$,
and  has the same distribution as $\X_{n-1}$ defined earlier.

To define $\H_n$, let $\Po^n$ be the restriction of $\Po$
to the set
 $$
\{ (x,t): t \leq  n f(X) \}.
$$
Let $\H_n$ 
be the  image of   the point set $\Po^n $
 under the mapping
$$
(x,t) \mapsto n^{1/d} (x - X).
$$
Given $X= x$, the point process $\Po^n $ is a
homogeneous Poisson process of intensity 1 on $\R^d \times [0, n
f(x)]$, and therefore, given $X= x$, $\H_n$ is a homogeneous
Poisson process on $\R^d$ of intensity $f(x)$. Define
 $\zeta_n $ 
to be the limit $\xi_\infty(\H_n)$.

Then $\H_n$ is a Cox process, where
 the randomness of the intensity measure comes from the
value of $f(X)$. Note  that the distribution of $\H_n$, and
hence that of $\zeta_n$, does not depend on $n$.

\begin{lemm}
\label{lemar27} Given  $K>0$, we have
 \bea \limn P[ n^{1/d}(\X'_{n-1}- X)
\cap B(\0;K) = \H_n\cap B(\0;K) ] = 1. \label{mar25a}
 \eea
\end{lemm}
{\em Proof.}
Suppose
 $X$ lies at a
Lebesgue point of $f$ (see e.g. \cite{Rudin}). Given $X=x$, the
expected  number of points of $\Po$ in $B(x; Kn^{-1/d}) \times
[0,\infty)$ that contribute to $\H_n$ but not to $\Po(n)$  is
$$
n\int_{B(x; Kn^{-1/d})} (f(x)-f(y))^+ dy
$$
which tends to zero because $x$ is a Lebesgue point
of $f$.  The expected
number of  points of $\Po$
in $B(x; Kn^{-1/d}) \times [0,\infty)$
that  contribute to $\Po(n)$ but not to  $\H_n$
is
$$
n\int_{B(x; Kn^{-1/d})} (f(y)-f(x))^+ dy
$$
which also tends to zero for the same reason.
Finally the probability that
$\Po(n) \cap B(x; Kn^{-1/d}) \neq \X'_{n-1} \cap
 B(x; Kn^{-1/d}) $
tends to zero as $n \to \infty$,
since
$|N(n)-(n-1)|$ is $o(n)$ in probability. Integrating over
possible values of $X$ and using Dominated Convergence, we obtain
(\ref{mar25a}).
$\qed$ \\

\begin{lemm}
\label{lemapr26} Suppose that $\xi$    is almost surely
stabilizing on $\Po_{\tau}$,  with limit
$\xi_\infty(\Po_{\tau})$, for all $\tau \in (0,\infty)$,
 and $\xi$ satisfies the  moments condition (\ref{bpm})
 for some $p > 1$.  Then
 \bea
 \label{apr26}
\limn     E[ n^{-1} H_{\xi_n}(\X_{n} ) ]
= \int_{\R^d} E [\xi_\infty (\Po_{f(x)})] f(x) dx. \eea
\end{lemm}
{\em Proof.}
Let $\eps >0$.  Then since
 $$ \xi_n(X; \X'_{n-1}) = \xi( \0; n^{1/d}(\X'_{n-1}- X)),
$$
for any $K\in \N$ we have \bea P[ |\xi_n(X; \X'_{n-1}) - \zeta_n|
> \eps] \nonumber
\\
\leq P[ n^{1/d}(\X'_{n-1}- X) \cap B(\0;K) \neq \H_n\cap B(\0;K)
] \label{m28a}
\\
+ P[ \overline{\xi}(\H_n;K) - \underline{\xi}(\H_n;K) > \eps].
\label{m28b} \eea
By the stabilization assumption,
we can choose $K>0$ so that (\ref{m28b}) is less than $\eps/2$,
and then by Lemma \ref{lemar27}, the expression
 (\ref{m28a})
 is also  less than $\eps/2$ for $n$ large.
Since $\zeta_n$ have the same distribution for all $n$, it follows
 that
\bea \xi_n (X; \X'_{n-1}) \tod
 \zeta_1.
\label{mar28} \eea Since the bounded $p$-th moments condition
(\ref{bpm}) is assumed to hold for some $p >1 $, the variables
$\xi_n (X; \X'_{n-1})$ are uniformly integrable, and hence their
expectations converge to that of
 $\zeta_1$.
By conditioning on $X$, we obtain \bea E[\zeta_1] =
 \int_{\R^d} E [\xi_\infty (\Po_{f(x)})] f(x) dx := \mu,
\label{mudef} \eea which is the right hand side of (\ref{apr26}).
Since $E[H_{\xi_n}(\X_n)]= n E[\xi_n(X; \X'_{n-1})]$,
this gives us (\ref{apr26}). $\qed$ \\

Next we consider $E[\xi_n(X_1;\X_n) \xi_n(X_2;\X_n)]$, and use a
refinement of the  coupling argument. Assume on a suitable
probability space that we have, independently, two
$d$-dimensional variables $X$ and $Y$ with density $f$, and two
homogeneous Poisson proceses $\Po $, $\Q$, both  of unit
intensity on $\R^d \times [0,\infty)$.

 Define coupled point processes
$\X'_{n-2}$ (a binomial process) and $\H_n^X$,  $\H_n^Y$ (both
Cox processes) and variables $\zeta_n^X$ and $\zeta_n^Y$,
all in terms of $\Po,\Q,X,$ and $Y$,
 as follows.
Let $\X'_{n-2}$ be obtained just as $\X'_{n-1}$ was  before, i.e.
 let $\Po(n)$ be the image of the  restriction of
$\P$ to the set $\{(x,t) \in \R^d \times [0\,\infty):
t \leq nf(x) \}$, under the projection $(x,t) \mapsto x$, and
let  $N(n)$ be the number of points of $\Po(n)$.
 Discard  $(N(n)-(n-2))^+$ of  the points of $\Po(n)$,
 chosen at random,
and add $(n-2-N(n))^+$ extra independent points with common density $f$.
The resulting set of points is denoted $\X'_{n-2}$
and has the same distribution as $\X_{n-2}$.

Let  $F_X$ be the half-space of points in $\R^d$ closer to $X$
than to $Y$, and let
  $F_Y$ be the half-space of points in $\R^d$
closer to $Y$ than to $X$. Construct $\H^X_n$ as follows. Let
$\Po^n_X$ 
 be  the restriction of $\Po$ to the set
$
F_X \times [0,n f(X)];
$
let $\Q^n_X $ be the restriction of $\Q$ to the set $F_Y \times
[0,nf(X)]$.
Let $\H_n^X$ be the  image of the  point process
 $ \Po^{n}_X  \cup \Q^n_X$
under the mapping
$$
(x,t) \mapsto n^{1/d} (x - X).
$$
Given $X= x$, the point process
 $ \Po^{n}_X  \cup \Q^n_X$
is a homogeneous Poisson process of intensity 1 on $\R^d \times
[0,nf(x)]$. Hence, given $X = x$,
 $\H^X_n$ is a homogeneous Poisson process
on $\R^d$ of intensity $f(x)$;
let $\zeta_n^X$ be the associated limit  $\xi_\infty(\H_n^X)$.

Construct $\H_n^Y$ in the following analogous manner.
 Let $\Po^{n}_Y$ be  the
restriction of $\Po$ to the set
$
F_Y \times [0,nf(Y)];
$
let $\Q_Y^n $ be the restriction of $\Q$ to the set $F_X
\times [0,nf(Y)]$. Let $\H_n^Y$ be the  image of the  point
process
 $ \Po^{n}_Y \cup   Q^n_Y$
under the mapping
$$
(x,t) \mapsto n^{1/d} (x - Y).
$$
By an argument similar to that used for $\H_n^X$, the point
process $ \H_n^Y$,  given $Y= y$, is a homogeneous Poisson
process on $\R^d$ of intensity $ f(y)$; we set $\zeta_n^Y : =\xi_\infty(\H_n^Y)$. \\

We can now prove the following result, which is the case $q=2$ of
Theorem \ref{generalLLN}.

\begin{prop}
\label{prop1}
 If $\xi$    is almost surely  stabilizing on $\Po_{\tau}$,  with limit
$\xi_\infty(\Po_{\tau})$, for all $\tau \in (0,\infty)$ and if
$\xi$ satisfies the  moments condition (\ref{bpm})
 for some $p > 2$,  then as $n \to \infty$,
 \bea
 \label{apr27}
  n^{-1 } H_{\xi_n}(\X_{n} )
\to \int_{\R^d} E [\xi_\infty (\Po_{f(x)})] f(x) dx \ \ \ {\rm
~in~}L^2.
\eea
\end{prop}
{\em Proof.} For any  $K>0$, we have \bea \limn P[
n^{1/d}((\X'_{n-2} \cup \{Y\})- X) \cap B(\0;K) = \H^X_n\cap
B(\0;K) ] = 1; \label{mar25a1} \eea \bea \limn P[
n^{1/d}((\X'_{n-2} \cup \{X\})- Y) \cap B(\0;K) = \H_n^Y
 \cap B(\0;K) ] = 1.
\label{mar25a2}
\eea
The proof of these facts
is just the same as that of  (\ref{mar25a}) using the
additional observation that
$$
\limn P[ B(X; Kn^{-1/d}) \subset F_X] = \limn P[ B(Y; Kn^{-1/d})
\subset F_Y]
 = 1.
$$

Note that $\H_n^X$ and $\H_n^Y$ are {\em independent} identically
distributed Cox processes. Independence follows by conditioning
on the values of $X, Y$; given these the point processes $\H^X_n$
and $\H_n^Y$ are constructed from Poisson processes on disjoint
regions of space. Therefore, for each $n$, the variables
$\zeta_n^X $ and $\zeta_n^Y$ are independent. Also the joint
distribution of $\zeta_n^X, \zeta_n^Y$ is independent of $n$.

Since by translation invariance
$$
\xi_n(X; \X'_{n-2} \cup \{X,Y\}) = \xi(\0; n^{1/d} ((\X'_{n-2}
\cup\{Y\})- X)),
$$
and
$$
\xi_n(Y; \X'_{n-2} \cup \{X,Y\}) = \xi( \0; n^{1/d} ((\X'_{n-2}
\cup\{X\})- Y)),
$$
it follows from  (\ref{mar25a1}) and  (\ref{mar25a2}), by a
similar argument to that which yielded (\ref{mar28}), that as
$n\to \infty$, \bea \xi_n (X; \X'_{n-2}\cup \{X,Y\})
 \xi_n (Y; \X'_{n-2}\cup \{X,Y\})
\tod \zeta_1^X \zeta_1^Y, \eea and hence \bea \xi_n (X_1; \X_{n})
\xi_n (X_2; \X_{n}) \tod \zeta_1^X \zeta_1^Y. \label{mar26b} \eea

By assumption,
 the bounded $p$-th moments condition
(\ref{bpm}) holds for some $p >2$. Then   by
Cauchy-Schwarz,
$$
\sup_{n \in \N} E[(\xi_n (X_1; \X_{n})
 \xi_n (X_2; \X_{n}) )^{p/2}] < \infty
$$
and therefore the variables $\xi_n (X_1; \X_n) \xi_n (X_2; \X_n)
$, defined for each $n \geq 2$, are uniformly integrable, so that
the convergence (\ref{mar26b}) also holds in the sense of
convergence  of means, i.e.
 \bea \limn E [ \xi_n (X_1; \X_{n})
\xi_n (X_2; \X_{n}) ] = E[ \zeta_1^X \zeta_1^Y] = \mu^2,
\label{mar27b} \eea
 with $\mu $ defined at (\ref{mudef}) Here we
have used independence of $\zeta_n^X $ and $\zeta_n^Y$.

To complete the proof, observe that 
\bean
 E[(n^{-1}H_{\xi_n}(\X_n)
)^2] = n^{-1}  E[\xi_n(X_1;\X_n)^2] 
+  (1 - \frac{1}{n}) E[\xi_n(X_1; \X_n) \xi_n(X_2; \X_n)],
\eean
and in the right hand side the first term tends to zero by the
bounded $p$-th moments condition, while  the second term tends
to  $\mu^2$ by (\ref{mar27b}). Therefore $E[n^{-1}H_{\xi_n}(\X_n)
] \to \mu$ and $ E[(n^{-1}H_{\xi_n}(\X_n) )^2] \to \mu^2, $ so
that $n^{-1}H_{\xi_n}(\X_n) $ tends to $\mu$
 in mean square. $\qed$ \\

Finally we  prove  $L^1$ convergence, i.e. the case $q=1$ of
 Theorem \ref{generalLLN}, completing the proof of that result.

\begin{prop}
\label{prop2}
 If $\xi$    is almost surely  stabilizing on $\Po_{\tau}$,  with limit
$\xi_\infty(\Po_{\tau})$, for all $\tau \in (0,\infty)$ and if
$\xi$ satisfies the  moments condition (\ref{bpm})
 for some $p > 1$,  then as $n \to \infty$,
 \bea
 \label{apr27a}
  n^{-1 } H_{\xi_n}(\X_{n} )
\to \int_{\R^d} E [\xi_\infty (\Po_{f(x)})] f(x) dx \ \ \ {\rm
~in~}L^1.
\eea
\end{prop}
{\em Proof.} Given $K>0$,  define the functional
$$
\xi^K(x,\X):= \min(\xi(x;\X),K).
$$
Then by the stabilization condition for $\xi$, the
truncated functional  $\xi^K$ also stabilizes on
 $\Po_\tau$ with limit
$\xi_\infty^K(\Po_\tau) := \min(\xi_\infty(\Po_\tau),K)$;
we leave the reader to verify this assertion.
Since $\xi^K$ is uniformly bounded, by Proposition
\ref{prop1} and the fact that $L^2$ convergence
implies $L^1$ convergence, we have
 \bea
 \label{apr27a2}
  n^{-1 } H_{\xi_n^K}(\X_{n} )
\to \int_{\R^d} E [\xi_\infty^K (\Po_{f(x)})] f(x) dx \ \ \ {\rm
~in~}L^1. \eea
 Moreover,
 \bean 0 \leq E[n^{-1}H_{\xi_n}(\X) -
n^{-1} H_{\xi_n^K} (\X)]
\\
= E[ \xi(n^{1/d}X_1; n^{1/d} \X_n) -
\xi^K(n^{1/d}X_1; n^{1/d} \X_n) ] \eean which tends to zero as $K
\to \infty$, uniformly in $n$, because the assumed moments
condition (\ref{bpm}), $p >1$, implies uniform integrability of
the family of variables $\xi(n^{1/d}X_1; n^{1/d} \X_n), $ $ n
\geq 1 $. Also, by monotone convergence the right hand side of
(\ref{apr27a2}) converges to the right side of (\ref{apr27a}) as
$K \to \infty$. Combining these facts,  and taking $K$ to
infinity in (\ref{apr27a2}), we obtain (\ref{apr27a}).$\qed$ \\

We now work towards a proof of
 Theorem \ref{Gtheo}. First we show that stabilization of
the graph $G$ as defined in Section \ref{secterm}
 implies an apparently stronger form of the stabilizing property.
\begin{lemm}
\label{stablem}
Suppose the graph $G(\X)$ is translation invariant and
stabilizes  on $\Po_1$.
Then with probability 1, for all points $X \in \Po_{1,0}$
there exists $R(X) < \infty$ such that the set of edges
of $G(\Po_{1,0})$
incident to $X$ is unchanged if points are added and/or deleted
outside $B(\0;R(X))$.
\end{lemm}
{\em Proof.} Let us say that a  point $x$ in a locally finite set
$\X$ is {\em unstable} for $G(\X)$  if there does {\em not} exist $r>0$ such that
$$
\Ed(x;G(\X ) ) = \Ed(x; G(\X \cap B(\0;r) )  \cup \A)
$$
for almost  all finite  $A \subset \R^d \setminus B(\0;r)$.

Let $t >1$.
By translation-invariance and the fact that
 the Poisson point process is its own
Palm distribution, the mean number of points of $\Po_1$ in
$B(\0;t)$  that are unstable for $G(\Po_1)$ is  equal to
$$
\int_{B(\0;t)} P[ x \mbox{ is unstable for } \Po_{1} \cup \{x\} ] dx
$$
which is zero by the assumption that $G$  stabilizes on $\Po_1$.
Therefore, since $t$ can be arbitrarily large, with probability 1
the Poisson process $\Po_1$ has no unstable points.

A further application of Palm theory for Poisson processes shows
that the mean number of pairs of distinct points $X,Y$ of $\Po_1$
in $B(\0;t)$ such that $Y$ is an unstable point of $\Po_1$, is
equal to
$$
\int_{B(\0;t)} \int_{B(\0;t)} P[ y    \mbox{ is unstable for } \Po_{1} \cup \{x,y\}]
dx dy
$$
and since the mean number of such pairs is zero,
 the above integral is zero, and therefore for
almost all $x$ in $B(\0;t)$,
\bea
\mbox{ Leb }\{y\in B(\0;t): P[ y
\mbox{ is unstable for }
\Po_1 \cup\{x,y\} \} ] >0\} =0.
\label{may13}
\eea
Choose $x_0$ in $B(\0;1)$, such that (\ref{may13}) holds
with $x=x_0$.
Then  the set of $z $ in $B(-x_0;t)$ such that
$$
P[z \mbox{ is unstable for } \Po_1 \cup \{0,z\} ] >0
$$
has zero measure.  Integrating over $z \in B(-x_0;t)$ and using
Palm theory yet again, we find that the mean number of Poisson
points  $X \in \Po_1$ in $B(-x_0;t)$ such that $X$ is unstable
for $G(\Po_1 \cup \{0\})$ is zero. Since  $t$ can be   arbitrarily
large, this
 gives us the result. $\qed$ \\

{\em Proof of Theorem \ref{Gtheo}.}
To prove (\ref{Gth1})  let $\xi(x;\X)$ be the reciprocal
of the order of the component containing $x$ in
$G(\X)$.  If  the component containing the origin
in $G(\Po_{1,0})$ is finite then by stabilization of
$G$ (Lemma \ref{stablem}) there  exists $R$ such that  alterations
to $\Po_1$ outside
$B(\0;R)$ will not cause any change in this component.
  If  the component containing the origin
in $G(\Po_{1,0})$ is infinite then given $\eps >0$ we can find a
connected subgraph of $G(\Po_{1,0})$ of order greater than
$\eps^{-1}$; then by Lemma \ref{stablem} there exists $R$ such
that  alterations to $\Po_1$ outside $B(\0;R)$ will not cause any
removal of edges  in this subgraph, and hence the functional
$\xi$ stabilizes on $\Po_1$ with limit $\sigma_G^{-1}$. By scale
invariance $\xi$  also stabilizes on $\Po_\tau$ with  limit having
the same mean as $\sigma_G^{-1}$. Since $\xi(x;\X)$ is uniformly
bounded by 1, the moments condition (\ref{bpm}) is trivially
satisfied and by Theorem \ref{generalLLN} we have (\ref{Gth1}).

To prove (\ref{Gth2}) let $\xi(x;\X)$ be equal to 1 if $G(\X)$
contains a subgraph isomorphic to $\Gamma$ with a vertex at $x$,
and to zero if not. This is bounded by 1, and by Lemma
\ref{stablem} $\xi$ stabilizes on  $\Po_1$ with limit equal to 1
if $E_G$ occurs and equal to 0 if not.
 Also $\xi$ stabilizes on $\Po_\tau$ with limit having the same
distribution, and Theorem \ref{generalLLN} applies to yield
(\ref{Gth2}).

To prove (\ref{Gth3}), observe that
 the  functional $L_\phi^G(\X) = \sum_{e \in G(\X)} \phi(|e|)$
 has the representation $L_\phi^G(\X) = L_{\xi}(\X) $ with
\bea
 \xi(x;\X) = \frac{1}{2} \sum_{e \in
 \Ed(x;G(\X))}
\phi(|e|).
\label{phiw2}
\eea
If $G$ is scale invariant,
then $G(\Po_\tau)$ has the same distribution as $G(\tau^{-1/d} \Po_1)$ and
therefore  $\xi$ stabilizes on $\Po_\tau$ with
\bea
E [ \xi_{\infty}({\cal P}_{\tau} )] = \frac{1}{2} E \sum_{e \in
\Ed(\0; G(\Po_{1,0}) )}
 \phi(\tau^{-1/d} |e|).
\label{psilindef}
\eea
Then (\ref{Gth3}) follows from Theorem \ref{generalLLN}.
$\qed$

%
%
\end{document}